# ON BEST PROXIMITY POINTS OF MULTIVALUED CYCLIC SELF-MAPPINGS ENDOWED WITH A PARTIAL ORDER

M. De la Sen

Institute of Research and Development of Processes. University of Basque Country

Campus of Leioa (Bizkaia) - Aptdo. 644- Bilbao, 48080- Bilbao. SPAIN

email: *manuel.delasen@ehu.es*

**Abstract**: The existence and uniqueness of fixed points of both the cyclic self-mapping and its associate composite self-mappings on each of the subsets are investigated if the subsets in the cyclic disposal are nonempty, bounded and of nonempty convex intersection is investigated.

**Keywords**: Best proximity point, cyclic self-mapping, fixed points, metric space, multivalued self-mapping, uniform convex Banach space.

## 1. Introduction

The extension of fixed point theory topics to the existence of either fixed points of multivalued self-mappings, [1-17], or common fixed points of several multivalued mappings or operators has received important attention. See, for instance, [13-17] and references therein. This paper investigates some properties of fixed point and best proximity point results for multivalued cyclic self- mappings under a general contractive-type condition based on the Hausdorff metric between subsets of a metric space [1], [4-6] and which includes a particular case the contractive condition for contractive single-valued self-mappings, [1-8] including the problems related to cyclic self-mappings. See, for instance, [4, 5, 9] and references there in. This includes strict contractive cyclic self -mappings and Meir-Keeler type cyclic contractions, [22-23]. Some fixed point results on contractive single and multivalued self-mappings, [1-2], [6-8], [18-19] and references therein, under some types of contractive conditions, have been revisited and extended in [1]. There is also a wide sample of fixed point type results available on fixed points and asymptotic properties of the iterations for self-mappings satisfying a number of contractive-type conditions while being endowed with partial order conditions. See, [16-17] and references therein. The objective of this research is the investigation of fixed point/ best proximity point results for multivalued cyclic self-mappings in complete metric spaces, or uniformly convex Banach spaces.

## 2. Some properties of multivalued cyclic self-mappings with a partial order

Assume that $(X,d)$ is a metric space for a set $X$ endowed with some metric $d:X\times X\to \boldsymbol{R}_{0+}$ with $\boldsymbol{R}_{0+}=\boldsymbol{R}_+\cup\{0\}$. Let $CL(X)$ be the family of all nonempty and closed subsets of the set $X$. If $A,B\in CL(X)$ then we can define $(CL(X),H)$ being the generalized hyperspace of $(X,d)$ equipped with the Hausdorff metric $H:CL(X)\to \boldsymbol{R}_{0+}$ induced by the metric $d:X\times X\to \boldsymbol{R}_{0+}$:

$$H(A,B)= max\left\{\sup_{x\in A} d(x,B),\ \sup_{y\in B} d(y,A)\right\} \qquad (2.1)$$

for two sets $A\subseteq X$ and $B\subseteq X$ which is finite if both sets are bounded and zero if they have the same closure. The distance between $A\subseteq X$ and $B\subseteq X$ is

$$D = d(A, B) = \inf_{x \in A, y \in B} d(x, y) = \inf_{x \in A} d(x, B) = \inf_{y \in B} d(y, A) \qquad (2.2)$$

Denote by $P(X)$, $B(X)$ and $CB(X)$ the sets of nonempty, nonempty and bounded and nonempty, and bounded and closed sets of $X$, respectively. The following relations hold:

$$D \leq H(A, B) \leq \delta(A, B) = \delta(B, A) = \sup_{x \in A, y \in B} d(a, b) \leq \delta(A, B) + \delta(B, C); \quad \forall A, B, C \in B(X) \qquad (2.3)$$

$$[(A, B \in CB(X)) \wedge H(A, B) < \varepsilon] \Rightarrow [\exists b \in B : d(a, b) < \varepsilon ; \forall a \in A] \qquad (2.4)$$

and $\delta(A, B) = 0$ if and only if $A = B = \{x\}$. Consider also a self-mapping $T : \bigcup_{i \in \bar{p}} A_i \to \bigcup_{i \in \bar{p}} A_i$, where $A_i$ are nonempty closed sets of $X$; $\forall i \in \bar{p} = \{1, 2, \ldots, p\}$, subject to the constraints $T(A_i) \subseteq A_{i+1}$ such that $A_{ip+j} \equiv A_j$ for any integer numbers $j \in [1, p-1] \cap \mathbf{Z}$ and $i \in \mathbf{Z}_{0+} = \mathbf{Z}_+ \cup \{0\}$ with $\mathbf{R}_{0+} = \mathbf{R}_+ \cup \{0\}$. If $p \geq 2$ then $T : \bigcup_{i \in \bar{p}} A_i \to \bigcup_{i \in \bar{p}} A_i$ is a $p$-cyclic self-mapping. If $p = 1$ then $T : A_1 \to A_1$ is, in particular, a self-mapping on $A_1$. We will also consider a partial order $\preceq$ on $X$ so that $(X, \preceq)$ is a partially ordered space and will assume, in general, that $T : \bigcup_{i \in \bar{p}} A_i \to \bigcup_{i \in \bar{p}} A_i$ is a multivalued $p$-cyclic self-mapping so that $A_i \ni x \to Tx (\neq \emptyset) \subset A_{i+1}; \forall i \in \bar{p}, \forall x \in \bigcup_{i \in \bar{p}} A_i$. The subsequent result does not assume a contractive condition for each iteration on adjacent subsets of the contractive mapping but a global contractive condition for the cyclic mapping for iterations on multiple strips of the $p$ subsets $A_i \subset X$; $i \in \bar{p}$. Therefore, the result that the distances between any two subsets being adjacent or not of [21] for nonexpansive self-mappings is not required. If $T : \bigcup_{i \in \bar{p}} A_i \to \bigcup_{i \in \bar{p}} A_i$ is a multivalued $p$-cyclic self-mapping then the set $BP(A_i) \subset A_i$ will be said to be the set of best proximity points between $A_i$ to $A_{i+1}$ if $d(A_i, A_{i+1}) = D_i = d(z, y)$ for all $z \in A_i$ and some $y \in Tz$. This concept generalized that of best proximity points in subsets of cyclic self-mappings established as follows. If $T : A_1 \cup A_2 \to A_1 \cup A_2$ is cyclic and single-valued then $x \in A_1$ and $Tx \in A_2$ are best proximity point if $d(A_1, A_2) = d(x, Tx)$, [20-21]. The following result extends a previous one for the case of non-cyclic self-mappings, [16-17]:

**Theorem 2.1**. Let $(X, \preceq)$ be a partially ordered space and $d : X \times X \to \mathbf{R}_{0+}$ with $(X, d)$ being a complete metric space. Let $A_i$ be a set of $p (\geq 2)$ nonempty, bounded and closed subsets of $X$ $X$; $\forall i \in \bar{p}$ (that is $A_i \in CB(X)$; $\forall i \in \bar{p}$) with $D_i = d(A_i, A_{i+1})$; $\forall i \in \bar{p}$ and let $T : \bigcup_{i \in \bar{p}} A_i \to \bigcup_{i \in \bar{p}} A_i$ be a multivalued $p$-cyclic self-mapping on $\bigcup_{i \in \bar{p}} A_i$ satisfying:

1. There exist $p$ real constants $k_i \in \mathbf{R}_{0+}$ satisfying $k = \prod_{i \in \bar{p}} [k_i] \in [0, 1)$ such that the following condition holds:

$$H(Tx, Ty) \leq k_i d(x, y) + (1 - k_i) D_i \qquad (2.5)$$

for any given $x \in A_i$ and $y \in A_{i+1}$ which fulfil $x \preceq y$, $\forall i \in \bar{p}$.

**2.** If $d(x,y) < d_0$ for some given $d_0 \in \mathbf{R}_+$, $y \in Tx$ and any given $x \in \bigcup_{i \in \bar{p}} A_i$ then $x \preceq y$ with $y \in A_{j+1}$ if $x \in A_j$ for any given $j \in \bar{p}$.

**3.** There are some $i \in \bar{p}$, some $x = x_i \in A_i$ and some $x_{i+1} \in Tx_i \subset A_{i+1}$ such that $d(x_i, x_{i+1}) < d_{0i}$ for some $d_{0i} > D_i$.

**4.** $d_0 \geq max\left(\max_{j \in \bar{p}} d_{0j}, \max_{j \in \bar{p}}\left(k_j(d_{0j} - D_j) + D_j\right)\right)$      (2.6)

Then, the following properties hold:

**(i)** There is a partially ordered subsequence $\hat{S}_i = \{x_{i+j+n_k p}\}_{n_k \in \mathbf{Z}_{0+}}$ of the partially ordered sequence $S(x_i) = \{x_{i+j}\}_{j \in \mathbf{Z}_{0+}}$, both of them of first element $x_i$, with respect to the partial order $(X, \preceq)$, such that $x_{i+j+n_k p} \in Q_{i+j}$ for $j \in \bar{p}$; $\forall k \geq k_0$, $n_k \in \mathbf{Z}_{0+}$ for some $k_0 \in \mathbf{Z}_{0+}$ and the given $i \in \bar{p}$, where $Q_{i+j} \subseteq Tx_{i+j-1} \subseteq T^{j-1}x_i \subseteq A_{i+j}$, for any $j \in \overline{p-1} \cup \{0\}$ and the given $i \in \bar{p}$, are $p$ closed "quasi-proximity" sets in-between each pair of adjacent subsets of the multivalued $p$-cyclic self-mapping $T : \bigcup_{i \in \bar{p}} A_i \to \bigcup_{i \in \bar{p}} A_i$ such that

$$D_{i+j} \leq d(x_{np+i+j+1}, x_{np+i+j}) \leq k_{i+j} D + (1 - k_{i+j})D = D \; ; \; \forall j \in \overline{p-1} \cup \{0\}, \forall n \in \mathbf{Z}_{0+} \quad (2.7)$$

where $D = \max_{j \in \bar{p}} D_j$ with $x_{np+i+j} \in Tx_{np+i+j-1} \subseteq A_{i+j}$; $\forall j \in \overline{p-1} \cup \{0\}$, $\forall n \in \mathbf{Z}_{0+}$ for the given $i \in \bar{p}$.

**(ii)** If $D_j = D$; $\forall j \in \bar{p}$ then any partially ordered sequence $S(x_i)$ of first element $x = x_i \in A_i$ fulfils:

$$\exists \lim_{n \to \infty} d(x_{np+i+j+1}, x_{np+i+j}) = D \quad (2.8)$$

; $\forall j \in \bar{p}$ and the given $i \in \bar{p}$, and $x_{np+j+1} \in Tx_{np+j} \subseteq A_{i+j+1}$; $\forall j \in \bar{p}$ (that is, $x_{np+j+1} \in A_{i+j+1}$ if $0 \leq j \leq p-i-1$ and $x_{np+j+1} \in A_{j-p+i+1}$ if $p-i < j \leq p-1$), $\forall n \in \mathbf{Z}_{0+}$. Let $BP(A_j)$ be the set of best proximity points between $A_j$ and $A_{j+1}$; $\forall j \in \bar{p}$. Then, there is a sequence $\{z_n^{(j)}\} \subset BP(A_j)$; $\forall j \in \bar{p}$ such that the following limit exists:

$$\lim_{n \to \infty} d(x_{np+j+1}, z_n^{(j)}) = D \; ; \; \forall j \in \bar{p} \text{ with } x_{np+j+1} \in Tx_{np+j} ; \; \forall n \in \mathbf{Z}_{0+} \quad (2.9)$$

**(iii)** If Assumption 3 is removed and (2.6) in Assumption 4 is replaced by the stronger condition

**5.** $d_0 > max\left(\max_{j \in \bar{p}}\left(D_j + diam(A_j)\right), \max_{j \in \bar{p}}\left(k_j(d_{0j} - D_j) + D_j\right)\right)$      (2.10)

then, Properties (i)-(ii) hold for any $x \in \bigcup_{i \in \bar{p}} A_i$.      □

Note that (2.5) is not guaranteed to be a cyclic contractive condition for each restricted map $T : \left(\bigcup_{j \in \bar{p}} A_j\right)|A_i \to \left(\bigcup_{j \in \bar{p}} A_j\right)|A_{i+1}$ since all the constants are not required to be less than one in (2.5) and, furthermore, (2.5) and Assumption 3 are fulfilled for some first element $x_i \in A_i$, $x_{i+1} \in Tx_i \subseteq A_{i+1}$ and some given $i \in \bar{p}$ in the partial order $(X, d)$. Note also that sequences fulfilling the partial order of Theorem 2.1 can always be built through iterations with the multivalued $p$-self-mapping for any arbitrarily chosen $A_i$ for any $i \in \overline{p}$ from (2.6) characterizing Assumption 4 of Theorem 2.1. Now, a particular case of Theorem 2.1 is stated:

**Theorem 2.2**. In addition to Assumptions 1-4 of Theorem 2.1, assume also:

6. $D_j = D = 0$; $\forall j \in \overline{p}$ (that is, $\bigcap_{j \in \overline{p}} A_j \neq \varnothing$)

7. The limit $x$ of any converging nondecreasing sequence $\{x_n\}_{n \in \mathbf{Z}_{0+}}$ is comparable to each $x_n$; $\forall n \in \mathbf{Z}_{0+}$ in the partial order $(X, \preceq)$, that is,

$$[x_n \preceq x(\neq x_n) \text{ for } x \in A_j, x_n \in A_j; \forall j \in \overline{p}, \forall n \in \mathbf{Z}_{0+}] \Rightarrow H(Tx, Tx_n) > k_i d(x, x_n) \quad (2.11)$$

Then, there is a sequence $\{x_{np+i+j}\}_{n \in \mathbf{Z}_{0+}}$ satisfying $x_{np+i+j} \in T^{np+j} x_i$ for some given initial element $x = x_i \in A_i$ and some given $i \in \overline{p}$; $\forall j \in \overline{p-1} \cup \{0\}$ which is non-decreasing and ordered with respect to the partial order $(X, \preceq)$ and fulfils the following properties:

(i) $\exists \lim_{n \to \infty} d(x_{np+i+j+2}, x_{np+i+j+1}) = 0$; $j \in \overline{p-1} \cup \{0\}$ and the given $i \in \overline{p}$ with $x_{np+j+2} \in T x_{np+j+1}$; $j \in \overline{p-1} \cup \{0\}$, $\forall n \in \mathbf{Z}_{0+}$ and the sequence $\{x_{np+i+j}\}_{n \in \mathbf{Z}_{0+}}$ is a Cauchy sequence; $j \in \overline{p-1} \cup \{0\}$.

(ii) The sequence $\{x_{np+i+j}\}_{n \in \mathbf{Z}_{0+}}$ for any $j \in \overline{p-1} \cup \{0\}$ and the given $i \in \overline{p}$ converges to a limit $\overline{x}$ in $\bigcap_{j \in \overline{p}} A_j$, which is a fixed point of the composite self-mapping $\hat{T}_j : A_j \to A_j$, where $\hat{T}_j = T^p = T \circ T \circ \cdots T$ ($p$ times) $= T^p|_{A_j}$ of domain $A_j$; $\forall j \in \overline{p}$ and also a fixed point of the self-mapping $T : \bigcup_{i \in \overline{p}} A_i \to \bigcup_{i \in \overline{p}} A_i$, that is, $\overline{x} \in \hat{T}_j \overline{x} \left(\subseteq \bigcap_{j \in \overline{p}} A_j\right)$ and $\overline{x} \in T^p \overline{x} \left(\subseteq \bigcap_{j \in \overline{p}} A_j\right)$; $\forall j \in \overline{p}$.

(iii) If, in addition, $(X, d)$ is a convex metric space, what holds, in particular, if $X$ is a Euclidean vector space and $d : X \times X \to \mathbf{R}_{0+}$ is the Euclidean metric, and $\bigcap_{j \in \overline{p}} A_j$ is convex, then $\overline{x} \in T \overline{x} \left(\subseteq \bigcap_{j \in \overline{p}} A_j\right)$ is the unique fixed point of $T : \bigcup_{i \in \overline{p}} A_i \to \bigcup_{i \in \overline{p}} A_i$ and $\hat{T}_j : A_j \to A_j$; $\forall j \in \overline{p}$ and also the unique fixed point of $T^p : \bigcup_{i \in \overline{p}} A_i \to \bigcup_{i \in \overline{p}} A_i$.

(iv) If Assumption 4 of Theorem 2.1 is replaced by Assumption 5 then Properties (i)-(iii) hold for any $x \in \bigcup_{i \in \overline{p}} A_i$.

(v) If $X$ is a Euclidean vector space then Property (iii) holds also if the condition of $(X, d)$ being a convex metric space is removed. □

## 3. The main result on best proximity points for non-intersecting subsets

An "ad hoc" version of Theorem 2.2 is being obtained in this section for the case of nonintersecting subsets by proving the convergence to unique best proximity points within each subset $A_i$, which are also $p$ respective unique fixed points of each of the composed self-mappings $\hat{T}_i : A_i \to A_i$ ; $\forall i \in \overline{p}$ if $(X, \|\ \|)$ is a uniformly convex Banach space endowed with the partial order $\preceq$ and the subsets $A_i$ ; $\forall i \in \overline{p}$ are nonempty, closed and convex sets.

**Theorem 3.1**. Let $T : \bigcup_{i \in \overline{p}} A_i \to \bigcup_{i \in \overline{p}} A_i$ be a multivalued $p(\geq 2)$-cyclic self-mapping on $\bigcup_{i \in \overline{p}} A_i$ with $A_i \in CB(X) \subseteq X$ ; $\forall i \in \overline{p}$ being all nonempty and convex with $D_i = d(A_i, A_{i+1})$ ; $\forall i \in \overline{p}$. Assume the following:

1. Let $X$ be a vector space and let $(X, d)$ be a convex complete metric space with $d : X \times X \to \mathbf{R}_{0+}$ being a homogeneous translation-invariant metric which induces a norm $\|\ \|$ on $X$ such that $(X, \|\ \|)$ is a Banach space.

2. $(X, \|\ \|)$ is a uniformly convex Banach space with metric convexity.

3. The complete metric space $(X, d)$, equivalently, the Banach space $(X, \|\ \|)$, is endowed with a partial order $\preceq$ defined by (2.5) with $x = x_i (\in A_i) \preceq y \in Tx (\subseteq A_{i+1})$ for any $(x, y) \in A_i \times A_{i+1}$ and some given $i \in \overline{p}$ such that the resulting $(X, \preceq)$ partially ordered space is subject to Assumptions 1-4 of Theorem 2.1 and Assumption 7 of Theorem 2.2.

Then, the following properties hold:

(i) There are unique best proximity points $\overline{x}_{j+1} \in T \overline{x}_j \subseteq A_{j+1}$ with $d(\overline{x}_j, \overline{x}_{j+1}) = d(\overline{x}_j, T \overline{x}_{j+1}) = D_j$, for each $j \in \overline{p}$ which are also unique fixed points of each of the restricted composite self-mappings $\hat{T}_j (\equiv T^p |A_j) : A_j \to A_j$ ; $\forall j \in \overline{p}$.

(ii) Take any $x = x_i (\in A_i) \preceq y = x_{i+1} \in Tx_i$ for any given $i \in \overline{p}$ (i.e. $x$ and $y$ are partially ordered with respect to the partial ordered set $(X, \preceq)$ and consider the partially ordered sequences $\{x_{np+j}\}$, being nondecreasing with respect to $\preceq$ while satisfying $x_{np+j+1} \in T x_{np+j}$ ; $\forall j \in \overline{p}$ of first element subject to $x = x_i (\in A_i) \preceq y = x_{i+1} \in Tx_i$ for any given $i \in \overline{p}$. Then, each of such sequences $\{x_{np+j}\}$ converges to the unique best proximity point $\overline{x}_j$ in $A_j$; $\forall j \in \overline{p}$ which is also the unique fixed point of each of the restricted composite self-mapping $\hat{T}_j : A_j \to A_j$. If $\bigcap_{i \in \overline{p}} A_i \neq \varnothing$ then $\overline{x} = \overline{x}_j \in \bigcap_{i \in \overline{p}} A_i$ is the unique fixed point of $T : \bigcup_{i \in \overline{p}} A_i \to \bigcup_{i \in \overline{p}} A_i$, $\hat{T}_j (\equiv T^p |A_j)$ and a fixed point of $T^p : \bigcup_{i \in \overline{p}} A_i \to \bigcup_{i \in \overline{p}} A_i$ ; $\forall j \in \overline{p}$.

(iii) If Assumption 4 of Theorem 2.1 is replaced by its Assumption 5 then the convergence to the above unique best proximity points holds for partially ordered sequences of first element $x \in \bigcup_{i \in \overline{p}} A_i$.

*Proof*: Note from the various hypothesis the uniformly convex Banach space $(X, \|\ \|)$ possesses the metric convexity property with respect to the norm metric $\|\ \|$ while it is endowed with a partial order $\preceq$ under Assumptions 1-4 of Theorem 2.1. From Property (ii) of Theorem 2.1, Eq. (2.8), the nonemptiness and closeness of the subsets $A_i \subseteq X$; $\forall i \in \overline{p}$ and Lemma 3.2 [(i)-(ii)], it follows that

$$\exists \lim_{n \to \infty} d(x_{np+i+j+1}, x_{np+i+j}) = D_{i+j}; \quad \exists \lim_{n \to \infty} d(x_{(n+1)p+i+j}, x_{np+i+j}) = 0; \quad \forall j \in \overline{p-1} \cup \{0\} \quad (3.1)$$

where $x_{np+i+j} \in Tx_{np+i+j-1} \subseteq A_{i+j}$, $x_{(n+1)p+i+j} \in T^p x_{np+i+j-1} \subseteq Tx_{(n+1)p+i+j-1} \subseteq A_{i+j}$; $\forall j \in \overline{p-1} \cup \{0\}$, $\forall n \in \mathbf{Z}_{0+}$ for the given $i \in \overline{p}$ and the iterated sequences; $\{x_{np+i+j}\}_{n \in \mathbf{Z}_{0+}}$; $\forall j \in \overline{p-1} \cup \{0\}$ and the given $i \in \overline{p}$ are partially ordered with respect to the partial order $\preceq$, from Theorem 2.1, of first element $x_{i+j} = x_j$ generated from the iteration $x_{np+i+j} \in Tx_{np+i+j-1}$; $\forall j \in \overline{p-1} \cup \{0\}$ and the given $i \in \overline{p}$ are all Cauchy sequences. Since $(X, d) \equiv (X, \|\ \|)$ is complete, it follows that $x_{np+i+j} \to \overline{x}_{i+j} \left( \in T^p x_{(n-1)p+i+j} \subseteq Tx_{np+i+j-1} \subseteq A_{i+j} \right)$ and $x_{np+i+j+1} \to \overline{x}_{i+j+1} \left( \in T^p x_{(n-1)p+i+j+1} \subseteq Tx_{np+i+j} \subseteq A_{i+j+1} \right)$ as $n \to \infty$; $\forall j \in \overline{p-1} \cup \{0\}$ and the given $i \in \overline{p}$ since $A_j \subseteq X$ is nonempty, bounded and closed; $\forall j \in \overline{p}$ and the given $i \in \overline{p}$. Thus, one gets from (3.1), since $A_j \subseteq X$ is nonempty, bounded and closed, and then boundedly compact, and also approximatively compact with respect to $A_{j-1}$ ([5],[27]), that:

$$D_{i+j} \leq d(x_{np+i+j+1}, x_{np+i+j}) \to d(\overline{x}_{i+j}, \overline{x}_{i+j+1}) = D_{i+j} = d(\overline{x}_{i+j}, T\overline{x}_{i+j}) \text{ as } n \to \infty$$

; $\forall j \in \overline{p-1} \cup \{0\}$ and the given $i \in \overline{p}$, where $\overline{x}_{i+j+1} \in T^p \overline{x}_{i+j+1-p} \subseteq T\overline{x}_{i+j}$; $\forall j \in \overline{p-1} \cup \{0\}$ and the given $i \in \overline{p}$. Since all the subsets $A_j \subset X$; $\forall j \in \overline{p}$ are nonempty, closed and boundedly compact; $\forall j \in \overline{p}$ then $\overline{x}_j \in A_j$ is a best proximity point in $A_j$ of $T : \bigcup_{i \in \overline{p}} A_i \to \bigcup_{i \in \overline{p}} A_i$ and it is also a fixed point of the restricted composite self-mapping $\hat{T}_j : \bigcup_{i \in \overline{p}} A_i \big|_{A_j} \to \bigcup_{i \in \overline{p}} A_i \big|_{A_j}$; $\forall j \in \overline{p}$. Thus, there are Cauchy, then convergent since $(X, d)$ is complete, sequences $\{x_{np+i+j}\}_{n \in \mathbf{Z}_{0+}}$ with respective first elements $x_{i+j} \in Tx_{i+j-1}$; $\forall j \in \overline{p-1} \cup \{0\}$ and the given $i \in \overline{p}$, each being convergent to $\overline{x}_{i+j} \in A_{i+j}$, such that $x = x_i$ is the first element of $\{x_{np+i}\}_{n \in \mathbf{Z}_{0+}} \subseteq A_i$ which consists of partially ordered elements with respect to the partial order $\preceq$ such that :

$$x_{i+j} \preceq \ldots \preceq x_{(n+1)p+i+j} \left( \in T^p x_{np+i+j} \subseteq Tx_{(n+1)p+i+j-1} \right) \preceq x_{(n+1)p+i+j} \preceq \ldots \preceq \overline{x}_{i+j} = \lim_{n \to \infty} x_{np+i+j} \quad (3.2)$$

with $\{x_{\ell p+i+j}\} \subseteq A_{i+j}$; $j \in \overline{p-1} \cup \{0\}$, $\forall \ell \in \mathbf{Z}_{0+}$ for the given $i \in \overline{p}$. But $\overline{x}_{i+j} \in A_{i+j}$; $\forall j \in \overline{p-1} \cup \{0\}$ and the given $i \in \overline{p}$, is a fixed point of the restricted composite self-mapping $\hat{T}_j : \bigcup_{i \in \overline{p}} A_i \big|_{A_j} \to \bigcup_{i \in \overline{p}} A_i \big|_{A_j}$; $\forall j \in \overline{p}$ and a fixed point of the composite self-mapping $T^p : \bigcup_{i \in \overline{p}} A_i \to \bigcup_{i \in \overline{p}} A_i$ to which the partially ordered sequences of first element $x = x_i \in A_i$ converge. It is also a best proximity point in $A_{i+j}$ of the self-mapping $T: \bigcup_{i \in \bar{p}} A_i \to \bigcup_{i \in \bar{p}} A_i$. The uniqueness property of each of those $p$ best proximity points $\bar{x}_j \in T\bar{x}_j$ in each of the subsets $A_j \subseteq X$ follows from their uniqueness as fixed points of the restricted self-mappings $\hat{T}_j : \bigcup_{i \in \bar{p}} A_i \big|_{A_j} \to \bigcup_{i \in \bar{p}} A_i \big|_{A_j}$ from Theorem 2.2 since $(X, d)$ is a convex metric space and the subsets $A_j \subseteq X$ are convex; $\forall j \in \bar{p}$. On the other hand, turns out that if all the subsets have nonempty intersection, such an intersection is convex so that the best proximity points are all identical and the unique fixed point of $T: \bigcup_{i \in \bar{p}} A_i \to \bigcup_{i \in \bar{p}} A_i$ and $T^p: \bigcup_{i \in \bar{p}} A_i \to \bigcup_{i \in \bar{p}} A_i$ from Theorem 2.2. This leads to the proofs of Properties (i)-(iii). □

**Remark 3.2**. (1) Theorem 3.1 proves the uniqueness of the best proximity points for any partially ordered sequences with first elements in any of the subsets of the multivalued $p$-cyclic self-mapping on $\bigcup_{i \in \bar{p}} A_i$ satisfying Assumptions 1-4 of Theorem 2.1 as it was commented in section 2 concerning such a theorem, the given $A_i \subset X$ for some $i \in \bar{p}$ to select the first two elements of the partial order can be chosen arbitrarily by construction. □

**Acknowledgements**
The author is very grateful to the Spanish Government by its support of this research trough Grant DPI2012-30651, and to the Basque Government by its support of this research trough Grants IT378-10 and SAIOTEK S-PE09UN12.